# Equitable Optimization of Patient Re-allocation and Temporary Facility Placement to Maximize Critical Care System Resilience in Disasters

Chia-Fu Liu[1]*, Ali Mostafavi[2]

[1] Ph.D Student, Zachry Department of Civil and Environmental Engineering, Texas A&M University, 199 Spence St., College Station, TX 77843-3136; e-mail: joeyliu0324@tamu.edu

[2] Associate Professor, Zachry Department of Civil and Environmental Engineering, Texas A&M University, 199 Spence St., College Station, TX 77843-3136; e-mail: amostafavi@civil.tamu.edu

## Abstract

End-stage renal disease patients face a complicated sociomedical situation and rely on various forms of infrastructure for life-sustaining treatment. Disruption of these infrastructures during disasters poses a major threat to their lives. To improve patient access to dialysis treatment, there is a need to assess the potential threat to critical care facilities from hazardous events. In this study, we propose optimization models to solve critical care system resilience problems including patient and medical resource allocation. We use human mobility data in the context of Harris County (Texas) to assess patient access to critical care facilities, dialysis centers in this study, under the simulated hazard impacts, and we propose models for patient re-allocation and temporary medical facility placement to improve critical care system resilience in an equitable manner.

The results show (1) the capability of the optimization model in efficient patient re-allocation to alleviate disrupted access to dialysis facilities; (2) the importance of large facilities in maintaining the functioning of the system. The critical care system, particularly the network of dialysis centers, is heavily reliant on a few larger facilities, making it susceptible to targeted disruption. (3) The consideration of equity in the optimization model formulation reduces access loss for vulnerable populations in the simulated scenarios. (4) The proposed temporary facilities placement could improve access for the vulnerable population, thereby improving the equity of access to critical care facilities in disaster. The proposed patient re-allocation model and temporary facilities placement can serve as a data-driven and analytic-based decision support tool for public health and emergency management plans to reduce the loss of access and disrupted access to critical care facilities and would reduce the dire social costs.

## I. Introduction

The objective of this study was to create an equitable optimization framework for patient re-allocation and temporary facility placement to maximize the resilience of critical care facilities network, with a focus on dialysis centers. Critical healthcare facilities like dialysis centers are crucial in safeguarding the wellbeing of patients with heightened vulnerability. The disruption of these services due to disasters can lead to perilous kidney failure in patients reliant on dialysis treatments (Smith et al., 2020). Patient risk is especially elevated during severe weather incidents, such as hurricanes, floods, or harsh cold conditions, when widespread kidney failure can result from interrupted access to these critical care facilities. Lempert & Kopp (2013) describe such a predicament as a "kidney failure disaster", an event that exposes a large number of patients, either on maintenance dialysis or recently diagnosed with acute kidney injury (AKI), at serious risk due to the unavailability of dialysis services. Historical data points to such health disasters, for instance, during Hurricane Katrina in 2005 (Bonomini et al., 2011; Kopp et al., 2007; Vanholder et al., 2009) and Hurricane Gustav in 2008 (Kleinpeter, 2009). For instance, the effects of Hurricane Sandy in 2012 were the major cause of kidney failure issues in the New York metropolitan area. Dialysis services were closed in anticipation of the storm or due to flooding, power outages, and structural damage caused by the storm (Lempert & Kopp, 2013). The closure of dialysis services in some severe flooded areas forced surrounding hospitals to house the evacuated patients, who were often admitted to emergency rooms with hyperkalemia. Despite the rapid response from renal communities, some patients are at increased health risk, and some may have suffered significant health consequences from missed dialysis sessions. The uncertainty and disruption caused by hazardous events may have resulted in acute and long-term mental health implications for maintenance dialysis patients. To build the resilience of the community of people with functional needs, it is important to establish a predetermined communication system to inform this population where they can receive dialysis treatment. Additionally, redundant communication methods and transport plans should be established. Notably, the disruption to routine dialysis sessions can have ripple effects, including an influx of patients to other dialysis centers, an increased strain on facilities caring for more dialysis-dependent patients, and more emergency department visits (Obialo et al., 2012; Saran et al., 2003).

---

* Corresponding author. Email: joeyliu0324@tamu.edu

The examination of disaster-induced disruption to vital dialysis centers remains an under-researched area within healthcare services and medical center studies. One of the rare investigations in this field, conducted by Kaiser et al. (2021), evaluated the flooding impact on dialysis centers in Harris County, Texas, during Hurricane Harvey, utilizing the flood maps from that weather incident. This study made use of flood zone categorizations provided by the Federal Emergency Management Agency (FEMA) to measure and classify dialysis centers based on their proximity to flood areas. However, focusing solely on the flood exposure of dialysis centers does not provide a comprehensive view of the potential threats to patients in the region arising from compromised access to these centers. Flooding can lead to multifaceted disruptions in accessing dialysis services, such as road inundation preventing patient travel (Redlener & Reilly, 2012); closures or malfunctions of dialysis centers due to facility flooding (Kaiser et al., 2021); and disturbances in the communities where dialysis-dependent patients reside (Du et al., 2012).

Two strategies in dealing with patients' disrupted access to critical care facilities, such as dialysis centers, include re-allocation of patients across the network of facilities in a region and setting up temporary facilities to meet the demand (Murakami et al., 2015). Different optimization methods have been proposed in the literature to solve patient and medical resource allocation problems (Fiedrich et al., 2000; Minciardi et al., 2009; Revelle & Snyder, 1995; Sun et al., 2014; Tsai et al., 2022; Ye et al., 2022; Yi & Özdamar, 2007). For the pandemic cases, Tsai et al. (2022) applied linear programming models to optimize the allocation of patients during the dengue fever epidemic. In the study, the objective function was to minimize the total travel distance of all patients. Sun et al. (2014) addressed patient and resource allocation between hospitals in a healthcare network during the pandemic influenza pandemic. The mathematical models take into account two objectives related to patients' cost of accessing healthcare services: (1) minimizing the total travel distance, and (2) minimizing the maximum distance a patient travels to a hospital. Ye et al. (2022) constructed a patient allocation model during major epidemics that considered the severity of patients' conditions by applying a multi-objective planning method.

For the disaster response cases, Minciardi et al. (2009) developed a mathematical model to assist decision makers in optimal resource allocation before and during a natural hazard emergency. Revelle & Snyder (1995) addressed emergency room location issues while respecting the maximum demand met. Fiedrich et al. (2000) investigated the allocation of available resources to the operational area to minimize the total death toll during the initial search and rescue phase after a major earthquake. Yi & Özdamar (2007) built an integrated location-distribution model to study the selection of temporary emergency centers that would result in maximum coverage of post-disaster medical needs in the affected area and optimal distribution of medical staff across both the temporary and permanent emergency response units.

Although past studies have implemented mathematical models in solving the problem of patient re-allocation and resource allocation, limited attention has been paid to the healthcare network optimization considering the possible infrastructure disruption in the aftermath of hazard events. Conversely, most studies make the assumption that existing facilities will not be affected by the disaster (Mete & Zabinsky, 2010; Mohammadi et al., 2016; Rabbani et al., 2016). This presumption, however, could be unrealistic, since the infrastructure, such as transportation facilities and medical facilities, could be severely damaged by a large-scale hazard event and remain inoperable for a period of time. Very few studies in the literature consider possible damage exclusively for the medical centers or the aid depots (Galindo & Batta, 2013; Huang et al., 2010; Paul & MacDonald, 2016).

Recognizing the gap, we propose a framework for disaster preparedness and response in healthcare networks considering infrastructure disruptions in the post-disaster period. Specifically, we focus on addressing the following research questions. (1) To what extent is the critical care facility network vulnerable to various infrastructure failure scenarios? (2) What is the optimized patient re-allocation plan for dialysis patients whose access is disrupted due to hazardous events? and (3) Where is a potential site for housing temporary medical facilities to improve access for socially vulnerable patients in an equitable manner? Accordingly, there are three objectives in the proposed model: the highest allocation effectiveness, the lowest transportation distance, and the equity of access to treatment for patients in each stricken area. The remainder of this paper is organized as follows. In the next section, we present the examined material and the formulation of the optimization models. In Section 3, numerical results of the studied case are presented to show how the model could help decision makers in determining patient allocation and the potential temporary facility placement in the healthcare system. In Section 4, the analysis and discussion based on the optimization results are presented. Section 5 contains concluding remarks. Fig. 1 presents the conceptual framework of this study.

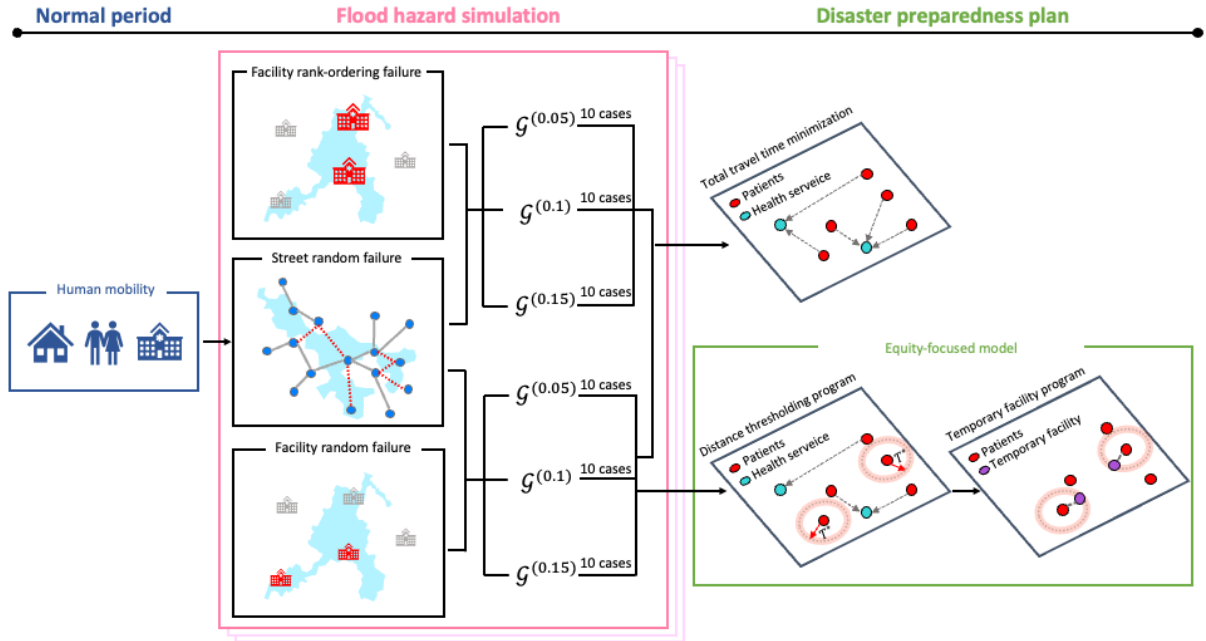

**Figure 1.** Conceptual framework of optimizing healthcare system resilience. The initial phase involves estimating the dialysis patient demand during the normal period (i.e., pre-disaster period). The subsequent flood-hazard simulation section examines 30 random failure scenarios which consist of facility random failure and street random failure, as well as 30 rank-ordering failure scenarios, which consist of facility rank-ordering failure and street random failure. Finally, the analysis includes total travel time minimization for both random failure and rank-ordering scenarios, while the equity-focused model is specifically applied to the random failure scenarios.

## II. Materials and Methods

### A. Population-facility Visitation Network and Demand Setting

This study uses the aggregated human mobility data to capture the dynamic visiting pattern of dialysis patients in the Houston metropolitan area. The human mobility dataset of stops at points-of-interest (dialysis centers in this study) from mobile devices, was collected from a mobility data provider. Each stopping point has been aggregated at the Census Block Group (CBG) level, thus forming the CBG-to-center visit. Dialysis demand exists in 2,010 CBGs out of a total of 2,144 CBGs in Harris County within which the Houston metro area is located. A total of 142 dialysis centers were included in the study. We used the two-week study period from August 1, 2017, to August 14, 2017, to estimate the number of patients in each CBG. A total of 5,308 visits were included in the two-week time window. These visits represent a sample of the actual number of visits. Since obtaining the actual number of visits is not feasible, we assume that these visits represent a fraction of the total visits. We discuss this assumption in the following section where we present the characteristics of facilities and their capacity.

### B. Topological datasets

To model the accessibility of the transport network to patients, we collect spatial data from OpenStreetMap, a collaborative mapping project that provides a free and publicly editable map. We imported the street network GIS data along with additional attributes (such as street type, street length, and speed limit) and then created the Harris County topological street network using the OSMnx package (Boeing, 2017). The road network, denoted by $\mathcal{G}$, characterizes intersections as nodes and road segments as edges in Harris County. To propose a proactive patient relocation plan considering the potential hazard event, we used the National Flood Hazard Layer (NFHL) to simulate flood hazards. As part of its National Flood Insurance Program, FEMA creates NFHL, consisting of digitized information for delineating floodplains in large geographic areas. The NFHL identifies not-at-risk zone as areas within the 500- or 100-year floodplains, as well as specially designated zones (e.g., coastal hazard zones). We extract the detailed floodplain boundary for our study case in Harris County.

## C. Flood Simulation

To simulate the impact of flooding on the Harris County dialysis healthcare network, we design two failure scenarios. We first identify the vulnerable zone based on the topological characteristics of the road network and the location of the medical facility. We identify the road segments overlaid with 100- and 500-year floodplains, denoted $S_r^{100}$ and $S_r^{500}$, respectively. In the same way, we identify the medical facility in the 100- and 500-year floodplain, denoted as $S_f^{100}$ and $S_f^{500}$. In the first scenario, the random failure scenario, we perform random removal of both road segments and medical facilities according to the flooding coefficient, $\delta$, defined as the percentage of flooded road segments in $S_r^{100}$ and the percentage of flooded facilities in $S_f^{100}$. For each setting of the flooding coefficient, $\mathcal{G}^{(\delta)}$, $\delta\%$ road segments in $S_r^{100}$ and $0.2 \times \delta\%$ road segments in $S_r^{500}$ would be randomly selected as flooded and thus inaccessible. In the same way, $\delta\%$ medical facilities in $S_f^{100}$ and $0.2 \times \delta\%$ medical facilities in $S_f^{500}$ are randomly selected as flooded and closed. In the second scenario, the capacity rank-ordering failure scenario, we perform the same random removal of the road segments but alternate the medical facility failure to capacity rank-ordering removal, where the medical facility on the floodplain, both $S_f^{100}$ and $S_f^{500}$, with the highest 10% capacity will be identified as flooded and therefore closed.

## D. Notations

All relevant notations used in the formulations are listed in Table 1.

**Table 1.** Notations used in the paper.

| General subscripts and sets | |
| --- | --- |
| $i, i'$ | Index of census block groups |
| $j$ | Index of medical centers |
| $l$ | Number of road segments |
| $m$ | Number of Census Block Groups |
| $n$ | Number of studied medical facility |
| $\mathcal{G}$ | Street topology network |
| $\mathcal{G}^{(\delta)}$ | Flooded street topology network with flooding coefficient $\delta$ |
| $S_r$ | Set of road segments, $S_r = \{1, \dots, l\}$ |
| $S_r^{100}, S_r^{500}$ | Set of road segments intersects with 100- and 500-year floodplains accordingly, $S_r^{100}, S_r^{500} \subseteq S_r$ |
| $S_c$ | Set of Census Block Groups, $S_c = \{1, \dots, m\}$ |
| $S'_c$ | Set of socio-vulnerable Census Block Groups, $S'_c \subseteq S_c$ |
| $S_f$ | Set of medical facilities, $S_f = \{1, \dots, n+1\}$ |
| $S_f^{100}, S_f^{500}$ | Set of dialysis cares intersect with 100- and 500-year floodplains accordingly, $S_f^{100}, S_f^{500} \subseteq S_f$ |
| **Parameters** | |
| $\delta$ | Flooding coefficient, $\delta \in [0,1]$ |
| $T_{ij}^f$ | Shortest travel time for trips from Census Block Group $i$ to facility $j$, $i \in S_c$, $j \in S_f$ |
| $T_{ii'}^c$ | Shortest travel time for trips from Census Block Group $i$ to Census Block Group $i'$, $i, i' \in S_c$ |
| $T^*$ | Threshold value for shortest travel time |
| $p_i$ | Lost-access patient in Census Block Group $i$, $i \in S_c$ |
| $p_i^{dt}$ | Lost-access patient after distance thresholding program in Census Block Group $i$, $i \in S_c$ |
| $c_j$ | Remaining capacity in medical facility $j$, $j \in S_f$ |
| $\rho$ | Median household income poverty line |
| **Variables** | |
| $x_{ij}$ | Relocated patient from Census Block Group $i$ to facility $j$, $i \in S_c, j \in S_f$ |
| $T$ | Total travel time for all lost-access patients |

## E. Optimization Formulation

### 1) Total travel time minimization model

In this study, we propose two optimization models to improve the dialysis healthcare network in the face of natural disasters. The first is the model of minimizing total travel time. The goal is to minimize the travel time $T$ for patients with lost access, which is expressed as follows

$$minT = min \sum_{i=1}^{m} \sum_{j=1}^{n+1} T_{ij}^f \cdot x_{ij} \qquad (1)$$

In flood scenarios, disruption to road segments and medical facilities will result in some patients losing access while others may retain accessibility without regard to travel time. We assume that patients who still have access will continue their treatment at the same facility, while patients who loses access will be transferred to the nearby facility based on the shortest travel time. We calculate the patient's shortest travel time to the medical facility, $T_{ij}^f$, based on the flooded street graph, $\mathcal{G}^{(\delta)}$, in each simulated scenario. The travel time for the entire route was calculated taking into account the free-flow travel speed and the travel length of each road segment. In addition, the model uses a dummy facility $n+1$ to receive the unsatisfied demand. Assigning patients to the dummy facility results in a prohibitively long travel time, $T_{i(n+1)}^f$, to the objective function. There are three constraints in the model of minimizing total travel time. First, the patient demand constraints for each CBG are

$$\sum_{j=1}^{n+1} x_{ij} = p_i \quad \forall i \in S_c \tag{2}$$

where, $x_{ij}$ is the relocated number from CBG $i$ to facility $j$ and $p_i$ is the number of lost-access patients in CBG $i$. Also, the number of relocated patients should not exceed the remaining capacity of facility $j$. The facility capacity constraints are

$$\sum_{i=1}^{n} x_{ij} \leq c_j \quad \forall j \in S_f \tag{3}$$

where, $c_j$ is the remaining capacity in medical facility $j$. Finally, the relocation number is a non-negative integer variable.

$$x_{ij} \in \mathbb{Z}^+ \quad \forall i \in S_c, \forall j \in S_f \tag{4}$$

*2) Equity-focused model*

The second model proposed in this study consists of two optimization programs, the distance thresholding program and the temporary facility program, which form an equity-focused model. The primary goal of the equity-focused model is to prioritize socially vulnerable (e.g., low-income and minority) patients. In large-scale disasters, the impact is particularly severe for minority communities, the elderly, the economically disadvantaged, and those with chronic illnesses. This demographic pattern is consistent with the population living with end-stage renal disease. Financial constraints and being in disaster-prone areas make economically disadvantaged people more vulnerable to disasters. Relocating dialysis treatment facilities places a significant burden on these patients, particularly those who rely heavily on public transport. To reduce this burden, when planning temporary post-disaster medical facilities, it is important to prioritize locations that are more accessible for these patients. This equity-focused model is geared towards setting up temporary medical facilities close to socio-vulnerable patients, thereby minimizing the need to travel long distances on public transport to travel to relocated facilities.

a. Distance thresholding program

We implement the equity-focused model for the random failure scenario. In the first part of the equity-focused model, we perform the distance thresholding program that has the same objective function in the total travel time minimization model as shown in Eq. (1). The program adds an additional constraint along with the Eqs. (2), (3), and (4) which is

$$x_{ij} = 0 \quad \forall T_{ij}^f > T^*, \forall i \in S'_c, \forall j \in S_f \backslash \{n+1\} \tag{5}$$

The additional constraint makes moving a socio-vulnerable patient to a facility with a travel time, $T_{ij}^f$, above the threshold, $T^*$, an infeasible solution. We have defined the socio-vulnerable population as patients residing on the CBG with a median household income below the poverty line, $\rho$. Therefore, following the relocation of the distance thresholding program, most non-vulnerable patients will be reassigned to the available medical facility, while the majority of socio-vulnerable patients who are far from the relocated facility will enter the second part of the equity-based model.

b. Temporary facility program

In the second part of the equity-focused model, select the locations for temporary medical facilities. We specify the possible locations for temporary facilities at the centroid of all CBGs. The program was designed as a multi-objective optimization problem. The first objective is to maximize the number of relocated patients, which is expressed as follows:

$$max \sum_{i \neq i'} x_{ii'} \tag{6}$$

The second objective function is to minimize the total travel time from the lost-access patients to the temporary facilities located in the centroid of CBGs. The objective is defined by:

$$minT = min \sum_{\substack{i=1 \\ i \neq i'}}^{m} \sum_{\substack{i'=1 \\ i \neq i'}}^{m} T_{ii'}^c \cdot x_{ii'} \tag{7}$$

where, $T_{ii'}^c$ is the shortest travel time for trips from Census Block Group $i$ to Census Block Group $i'$. For each CBG, the relocated number should not exceed the remaining patients with lost access after the distance thresholding program indicated by $p_i^{dt}$. The demand constraint is shown in Eq. (8).

$$\sum_{j=1}^{m} x_{ii'} \leq p_i^{dt} \quad \forall i \in S_c \tag{8}$$

Also, the non-negative integer constraint still applies to the relocating variables.

$$x_{ii'} \in \mathbb{Z}^+ \quad \forall i, i' \in S_c \tag{9}$$

In this study, we solve the optimization problems using the PuLP package with the COIN Branch-and-Cut (CBC) solver in Python.

## III. Results

### A. Impacts of simulated flood hazard

To build the Harris County dialysis facility network and the demand for facilities, we first estimate the capacity of each dialysis center based on the assumption that the two-week demand accounts for 90% of each facility's total capacity. This assumption is mainly due to the infeasibility of obtaining actual demand and capacity data. The distribution of the estimated capacity for the dialysis center is shown in Fig. 2. The capacities of the dialysis centers follow a long-tail distribution with a small number of facilities having the largest capacity and greatest demand.

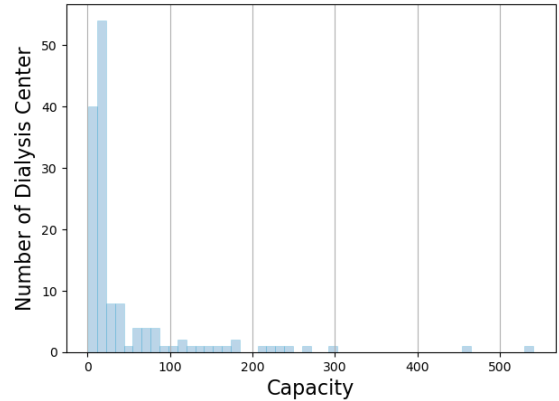

**Figure 2.** The distribution of estimated service capacity for the 142 dialysis centers in Harris County, Texas.

Second, we perform geospatial processing on the road network and the floodplain to classify the road segment. Of the 354,546 total road segments in Harris County, Texas, 54,084 segments denoted as $S_r^{100}$, that intersect with the 100-year floodplain, while 63,170 segments denoted as $S_r^{500}$ intersect with the 500-year floodplain. The geographic topology of $S_r^{100}$ and $S_r^{500}$ is shown in Fig. 3. In addition, of the 142 total Harris County medical facilities examined, there are 25 medical facilities, $S_f^{100}$, at the 100-year floodplain and 27 facilities, $S_f^{500}$, at the 500-year floodplain.

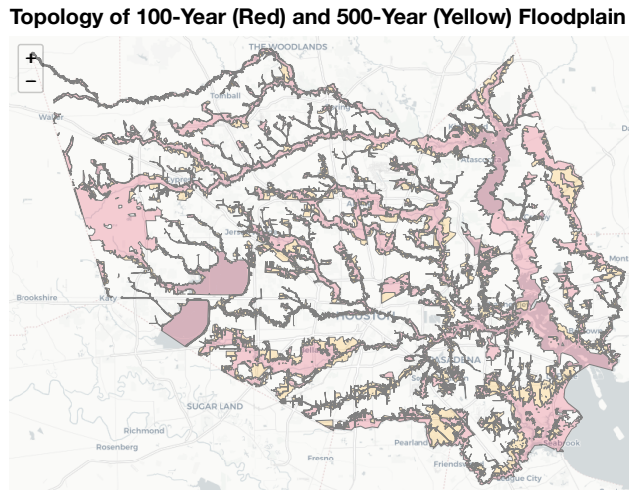

(a)

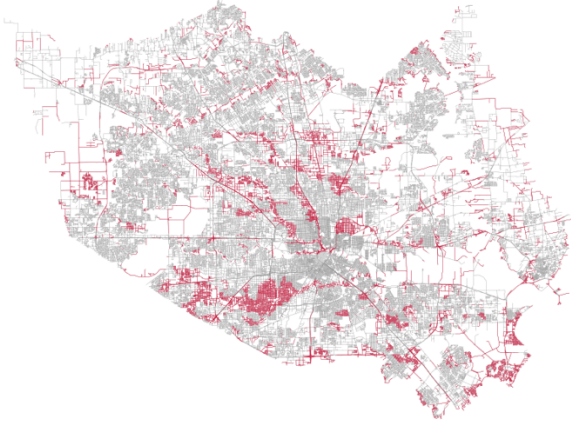
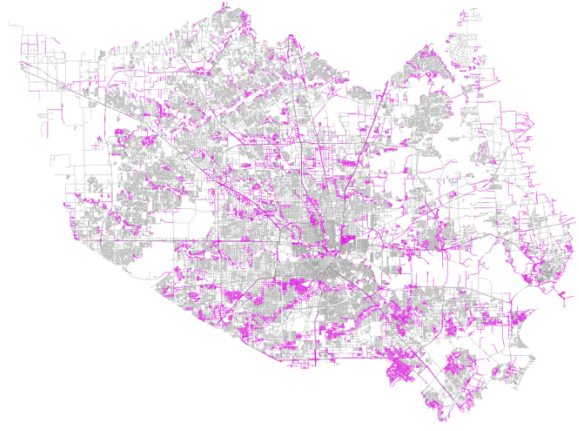

Road Segments intersect with 100-year Floodplain    Road Segments intersect with 500-year Floodplain

(b)                                                 (c)

**Figure 3.** (a) 100-year and 500-year floodplains in Harris County, Texas. (b) The geographic topology of the road segments intersects with the 100-year floodplain. (c) The geographic topology of the road segments intersects with the 500-year floodplain.

In the random failure scenario, we generated flooded street topology networks $\mathcal{G}^{(\delta)}$ by setting $\delta$ equals 0.05, 0.1, and 0.15 to simulate both road failure and medical facility failure. In this scenario, $\delta$% of $S_r^{100}$ and $S_f^{100}$ are randomly selected as flooded and $0.2 \times \delta$% of $S_r^{500}$ and $S_f^{500}$ are randomly selected as flooded. We simulated 10 cases for each $\delta$ setting. The average flooded road segment and medical facility in each $\delta$ setting is shown in Table 2.

| Flood coefficient | $\delta = 0.05$ | $\delta = 0.1$ | $\delta = 0.15$ |
|---|---|---|---|
| Road segment | 3,328.6 | 6,644 | 9,942.8 |
| Medical center | 1 | 3 | 5 |

**Table 2.** The average flooded road segment and flooded medical facility in random failure scenario under different $\delta$ settings.

In the capacity rank-ordering failure scenario, the identification of the flooded road segment follows the same procedure as in the random failure scenario, resulting in an identical $\mathcal{G}^{(\delta)}$. However, for the medical facility failure, we select facilities in the floodplain with the largest 10% capacity, which is 5 of 52 facilities to be flooded. The goal of the capacity rank-ordering failure scenario is to stress test the system and assess the level of dependency of the Harris County dialysis community on these major dialysis centers. Similarly, to assess the potential impact on the dialysis healthcare system in Harris County, we generate ten cases for each flooded street network $\mathcal{G}^{(\delta)}$, giving a total of 60 flooding cases.

We run the flooding simulations; the results show that in the random failure scenario, the average access-lost patient under three flooding coefficient settings is 95 ($\mathcal{G}^{(0.05)}$), 324 ($\mathcal{G}^{(0.1)}$), and 480 ($\mathcal{G}^{(0.15)}$). Meanwhile, in the capacity rank-ordering failure scenario, the average loss of access for patients under three flooding coefficient settings is 1,342 ($\mathcal{G}^{(0.05)}$), 1396 ($\mathcal{G}^{(0.01)}$), and 1,483 ($\mathcal{G}^{(0.15)}$). The distributions of lost-access distributions are shown in Fig. 4. Lost access means the patients would not be able to access any facility in the region since all facilities are out of capacity or out of service. Disrupted access, on the other hand, means patients need to take longer travel to access facilities.

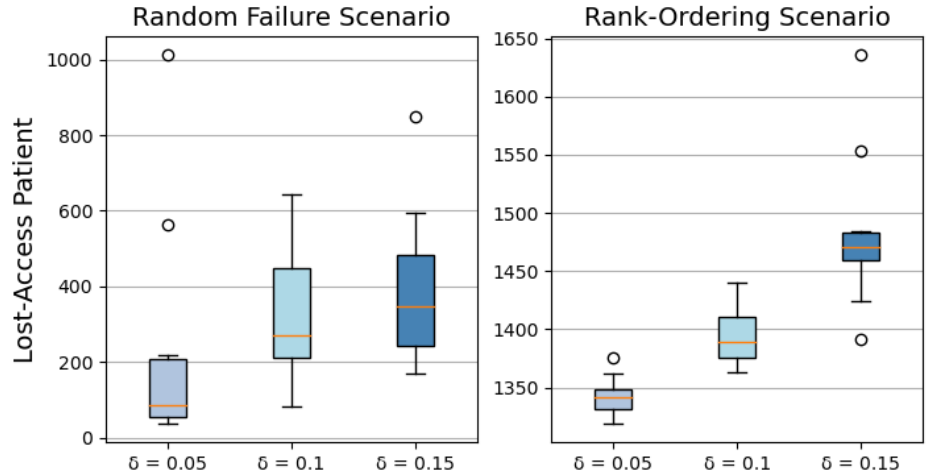

**Figure 4.** The distribution of patients losing access to facilities in random failure and capacity rank-ordering failure scenario.

In the random failure scenario, the average travel time (in seconds) of patients still able to access their medical facility under the impact of flooding is 435.27 ($\mathcal{G}^{(0.05)}$), 451.35 ($\mathcal{G}^{(0.1)}$), and 471.26 ($\mathcal{G}^{(0.15)}$). Meanwhile, in the capacity rank-ordering failure scenario, the average travel time of patients who still have access to their medical facility despite the flood impact is 443.26 ($\mathcal{G}^{(0.05)}$), 459.53 ($\mathcal{G}^{(0.1)}$), and 476.5 ($\mathcal{G}^{(0.15)}$). The distribution of the average travel time at different flooding coefficients in each flood scenario is shown in Fig. 5.

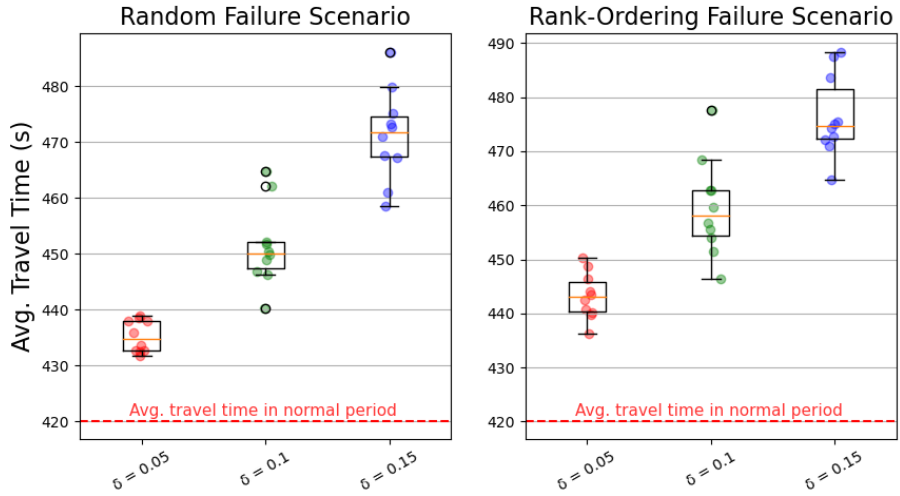

**Figure 5.** The distribution of average travel time with three flooding coefficients in random failure scenario and capacity rank-ordering capacity. The red dashed line with a value of 420.19 (seconds) indicates the average travel time in the normal period.

## B. Total travel time minimization model

In Fig. 6, we present the optimization result for minimizing the total travel time under random failure and facility capacity rank-ordering scenarios.

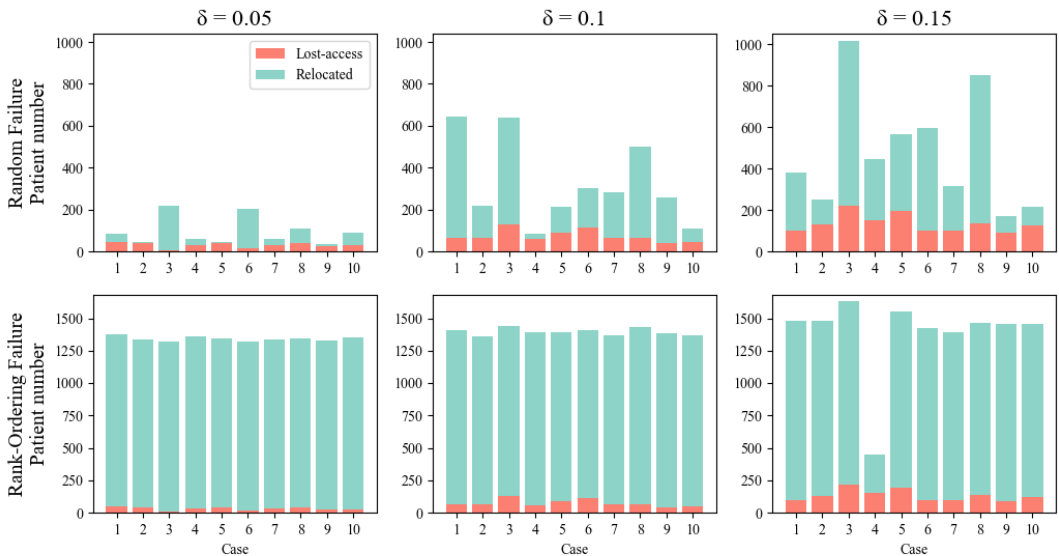

**Figure 6.** The optimization result of the total travel time minimization model under random failure scenario and capacity rank-ordering scenarios. The green area represents the patients who were re-allocated to nearby dialysis centers. The red area represents patients who still did not have access according to the total travel time minimization model.

## C. Equity-focused model

In this study, we set the poverty line, $\rho$, as the first quartile of the median household income of all patients, which is \$33,956.75. Patients in the CBG with a median household income of less than $\rho$ are classified as a socio-vulnerable population. We also set the travel time threshold, $T^*$, for the socio-vulnerable population as the median travel time of all patients in the normal period, which is 251.4 seconds. The relocation of socio-vulnerable patients with a travel time greater than $T^*$ is identified as an infeasible solution under the distance thresholding program. Fig. 7 shows the distribution of median house income for CBGs in Harris County and the distribution of travel time for all patients over the normal period.

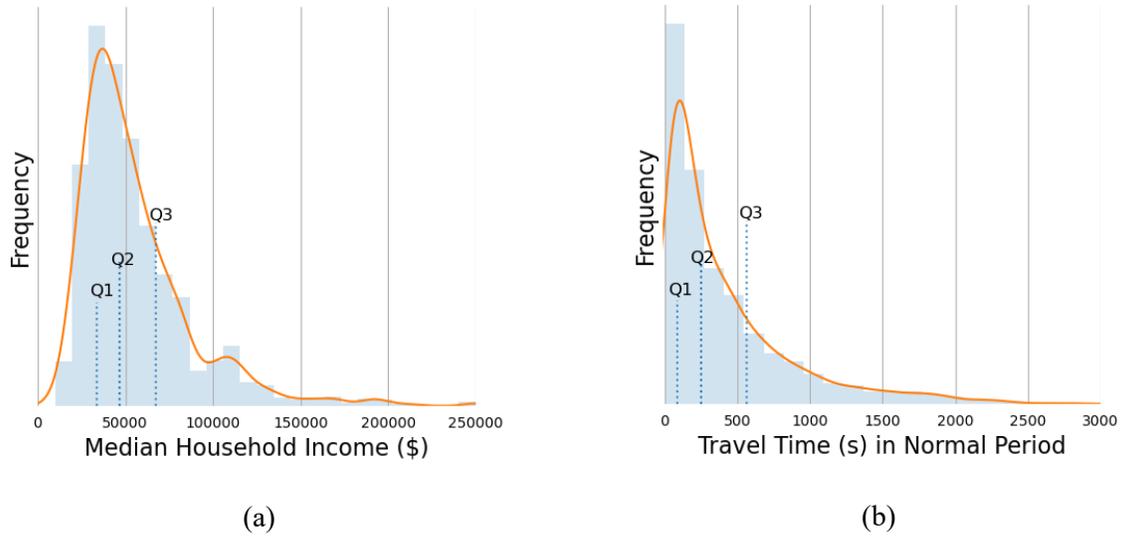

(a)                                                           (b)

**Figure 7.** (a) The distribution of dialysis patients' median household income distribution. (b) The distribution of travel time of all patients during the normal period.

We aggregate the ten cases in each flooding coefficient setting to represent the points of need for temporary facilities that could provide the shortest travel time for the nearby patients losing access. In the flooded road topology network $\mathcal{G}^{(0.05)}$, there are 18 temporary facility demand points with a total demand of 59. In addition, in the flooded street topology network $\mathcal{G}^{(0.1)}$, there are 56 temporary facility demand points with a total demand of 183. Finally, in flooded street topology network $\mathcal{G}^{(0.15)}$, there are 80 temporary facility demand points with a total demand of 281. Fig. 8 shows the geographic distribution of the aggregated temporary facility demand point with ten cases under different flooding coefficients.

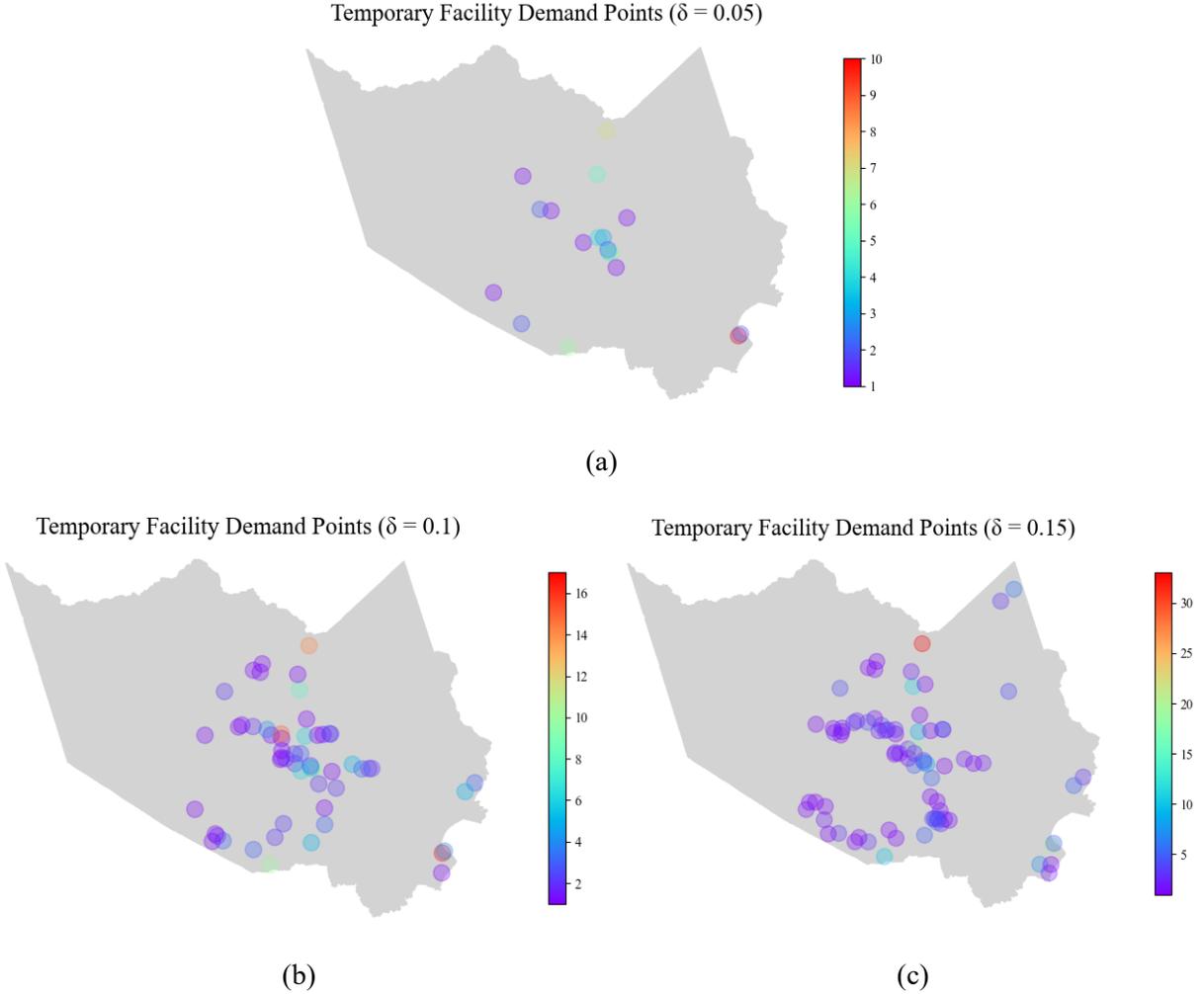

**Figure 8.** The geographic distribution of demand points for temporary medical facilities under three flooding coefficients with their ten-case aggregate demand.

## IV. Analysis and Discussion

In this study, we designed two failure scenarios to assess the extent to which the examined healthcare network is dependent on large medical facilities. As shown in Fig. 2, there are only 18 medical facilities with an estimated capacity greater than 100, which means that 87.32% of medical care has a capacity less than 100. The result shows that the dialysis healthcare network relies heavily on a few medical centers with large capacities. As shown in Fig. 9, the critical care facility capacity rank-ordering results in a higher average travel time (green) than in the random failure scenario (red). This shows that the accessibility of dialysis is highly dependent on the operation of these large medical centers. This means that once these large medical facilities are perturbed, patients will find it difficult to find an alternate facility nearby to continue their maintenance dialysis. In other words, the critical care facility network has a scale-free structure and is vulnerable to targeted attacks on hub nodes (i.e., large facilities).

Although the capacity rank-ordering medical facility failure results in many more patients losing access (red plus green area) as shown in the second row of Fig. 6, the total travel time minimization model could re-allocate most of the patients (green) and achieve the same level of performance, in terms of the number of remaining patients with lost access (red), as in the random failure scenario.

In addition, as shown in Fig. 9, the model for minimizing total travel time will re-allocate the socio-vulnerable patient with a longer average travel time, particularly with less variance in the rank-ordering failure scenario (green) than in the random failure scenario (red). By formulating the equity-focused model, we can greatly reduce the average travel time (blue) of socially vulnerable patients in all flooding coefficient settings, thus improving the equity of accessibility of dialysis patients when the facility network is perturbed.

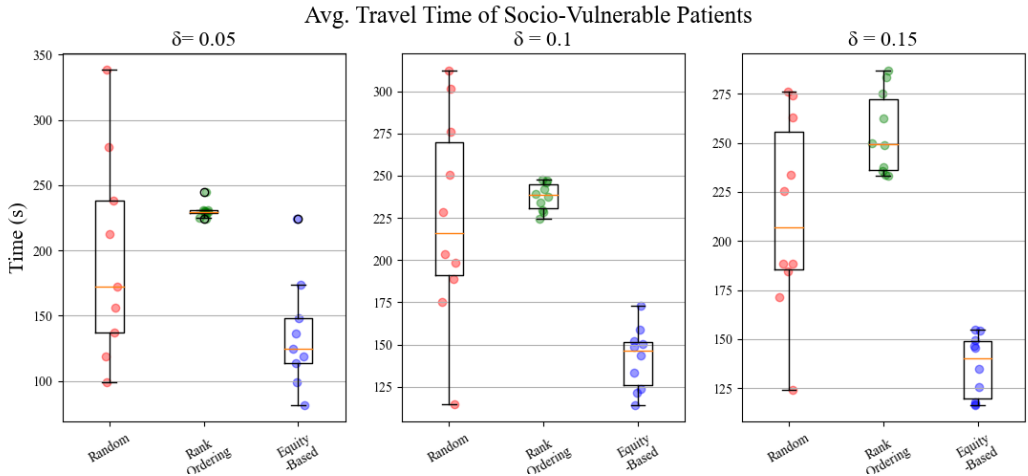

**Figure 9.** The average travel time with the travel minimization model (red and green) and the equity-based model (blue) under three flooding coefficients.

From Fig. 8 we can observe that in the central part of Harris County, there is a strip of demand points in all three flooding coefficient settings, regardless of its demand magnitude. This could be an indicator of high patient demand and an alarm signal that the central-area road segment is highly vulnerable to simulated flooding. In addition, although there are some large demand points in the northern periphery, particularly under $\mathcal{G}^{(0.15)}$ setting, it could be a false positive signal, in particular considering we have excluded the possible medical center north of the border of Harris County. The boundary of the study region is the limits of Harris County. Patients living in the boundary regions of the county might visit facilities in the neighboring county. Hence, demand points identified in the periphery of the county should be further examined in light of proximity of facilities in the neighboring county. Also, the facilities in neighboring counties are likely to be impacted by the flood event as well. Hence, not considering the neighboring county's facilities does not undermine the results of the optimization model.

## V. Concluding Remarks

This study proposed an equity-focused optimization framework to assess the critical care facility network resilience in disasters. The study and its outcomes have multiple important contributions. First, the proposed optimization framework is among the first efforts to characterize and improve the resilience of regional critical care facility networks in disasters. The framework proposed in this study can serve as a decision support tool to inform emergency managers and public health officials to better understand, prepare, and respond to the effects of disasters on dialysis centers by optimizing patient re-allocation and temporary facility placement. Second, this study incorporates equity in the optimization model formulation which is mostly ignored in prior studies. Disasters are known to disproportionately impact vulnerable populations; if decision support models (such as the optimization model presented in this paper) do not consider equity aspects, the impacts on vulnerable populations will be exacerbated. Third, this study developed the optimization model of population-facility network based on observational location-based data that would provide a more realistic representation of patients' dependence on different facilities.

Specifically, we simulated the disruption of the road network and the closure of dialysis facilities based on different levels of flood severity. The results show that: (1) the critical care facility network is highly dependent on certain large dialysis centers, which means that the system has the characteristics of scale-free networks and is vulnerable to targeted disruption such as capacity rank-ordering failure. (2) In addition, by assessing the geographic distribution of temporary facility demand points, we also identified the areas of the dialysis patient community that are vulnerable to critical facility failures. A possible solution is to develop a distributed dialysis healthcare system in the study area. While a centralized healthcare system can leverage economies of scale to provide healthcare services more cost-efficiently, especially in areas with large populations, a distributed healthcare system could allow better access to patients. By providing dialysis care closer to where patients live, a distributed healthcare system can reduce travel time and costs and make healthcare more accessible. In addition, a distributed critical care system can be more flexible and responsive to local needs and conditions, enabling providers to adapt to changing patient demands when facing hazardous events. Still, distributed facilities can be more complex to manage and may require greater investment in infrastructure. (3) Furthermore, this study identified potential sites for temporary medical facilities that could be deployed to enhance the healthcare system in the context of disaster resilience, with a particular focus on the socioeconomically vulnerable population. Given socio-economically vulnerable patients' transportation barriers to reach their relocated dialysis treatment, the increased travel time caused by relocation could pose a significant burden. Therefore, to promote an equitable healthcare network, when designing temporary medical facilities, the dialysis community should consider prioritizing locations that are more accessible to these socioeconomically vulnerable patients.

The use of data-driven methods for disaster risk management can significantly reduce the impacts of hazards on communities (and vulnerable populations in particular). This study is a first of its kind to leverage location-based datasets along with datasets related to critical care facilities, road networks, and hazard exposures in creating and testing an optimization model to efficiently re-distribute dialysis-dependent patients when access to a number of facilities is disrupted during extreme weather events. The current approach to preparing and responding to disruptions in access to dialysis facilities is rather reactive and could lead to sub-optimal allocation of patients to other facilities. Through the use of the proposed data-driven optimization model, public health officials and emergency managers could proactively evaluate different scenarios of road inundations and facility outages to identify: (1) areas most vulnerable to loss of access to dialysis centers; (2) optimal ways to allocate patients to dialysis centers for different scenarios individually and collectively; (3) capacity utilization of dialysis centers and facilities that play a critical role in the overall redundancy of the network of facilities within a city or region. These insights could inform plans to reduce the impacts of disrupted access by increasing the capacity of most important facilities during extreme weather events, building new facilities in areas with the greatest vulnerability, and reducing vulnerability of the existing facilities through backup power and water supply. These proactive measures obtained based on data-driven methods could prevent catastrophic kidney failure disasters in the face of extreme weather events. Also, the method and data presented in this study could be used for optimizing access to other critical care facilities beyond dialysis centers.

**Data availability**

The data used in this study are not publicly available under the legal restrictions of the data provider. Interested readers can request it from the provider directly.

**Code availability**

The code that supports the findings of this study is available from the corresponding author upon request.

**Acknowledgments**

This material is based in part upon work supported by the National Science Foundation under Grant CMMI-1846069 (CAREER). Any opinions, findings, conclusions, or recommendations expressed in this material are those of the authors and do not necessarily reflect the views of the National Science Foundation.

**Author Contributions**

C.F.L. and A.M. conceived the idea. C.F.L collected the data and carried out the analyses. C.F.L. and A.M. wrote the manuscript.